\documentclass{SCIYA2017enOL}
\online

\def\span{\hbox{\rm span$\,$}}

\begin{document}

\ensubject{fdsfd}

\ArticleType{ARTICLES}
\Year{2019}
\Month{}%
\Vol{60}
\No{}
\BeginPage{1} %
\DOI{}
\ReceiveDate{January 31, 2019}
\AcceptDate{May 5, 2019}

\title[]{Infinite dimensional Cauchy-Kowalevski and Holmgren type theorems}
{Infinite dimensional Cauchy-Kowalevski and Holmgren type theorems}

\author[1]{Jiayang YU}{jiayangyu@scu.edu.cn}
\author[2,$\ast$]{Xu ZHANG }{zhang$\_$xu@scu.edu.cn}


\AuthorMark{Jiayang Yu}

\AuthorCitation{Jiayang Yu, Xu Zhang}

\address[1]{ School of Mathematics, Sichuan University, Chengdu 610064, P. R. China}
\address[2]{ School of Mathematics, Sichuan University, Chengdu 610064, P. R. China}


\abstract{The aim of this paper is to show Cauchy-Kowalevski and Holmgren type theorems with infinite number of variables. We adopt von Koch and Hilbert's definition of analyticity of functions as monomial expansions. Our Cauchy-Kowalevski type theorem is derived by modifying the classical method of majorants. Based on this result, by employing some tools from abstract Wiener spaces, we establish our Holmgren type theorem.}



\maketitle

\section{Introduction}
The classical Cauchy-Kowalevski theorem asserts
the local  existence and uniqueness of analytic solutions to quite general partial differential equations with analytic coefficients
and initial data in the finite dimensional Euclidean space $\mathbb{R}^{n}$ (for any given $n\in \mathbb{N}$).
In 1842, A. L. Cauchy first proved this theorem for the second order case; while in 1875, S. Kowalevski proved the general result. Both of them used the method of majorants.
Because of its fundamental importance, there exists continued interest to generalize and/or improve this theorem (See \cite{Cha, Fri, Lascar79, Mog, Nir, Nis, Ovs, Saf, Tre, Yam92, Zub01} and the references cited therin).
In particular, some mathematicians studied abstract forms of Cauchy-Kowalevski theorem  in the context of Banach spaces and Fr\'{e}chet derivatives. For example, the linear cases were considered by L.~Ovsjannikov (\cite{Ovs}) and F.~Tr\`{e}ves (\cite{Tre}). Later, L.~Nirenberg (\cite{Nir}) obtained a nonlinear form, while T.~Nishida (\cite{Nis}) simplified Nirenberg's proof and obtained a more general version, M.~Safonov (\cite{Saf}) gave another proof of Nishida's theorem. There are also some abstract Cauchy-Kowalevski theorems in this respect (\cite{Cha, Mog, Yam92}).

On the other hand, Holmgren's uniqueness theorem states the uniqueness of solutions to linear partial differential equations with analytic coefficients in much larger class of functions than analytic ones. In 1901, E.~Holmgren (\cite{Hol}) first proved this theorem for the case of $n=2$ in the context of classic solutions, while in 1949, F.~John (\cite{Joh}) extended it to the general case of $n$ variables. Later the result was extended first to the setting of distribution solutions by L.~H\"{o}rmander (\cite{Hor1, Hor2}) and then to that of hyperfunction solutions by H.~Hedenmalm (\cite{Hed}). All of the above works are addressed to the case of finite dimensional spaces. In the 1970s, B.~Lascar (\cite{Lascar76, Lascar79}) gave a Banach space version of Holmgren's uniqueness theorem. Unfortunately, the proof of Lascar's result in this respect is incomplete and questionable, and therefore, about thirty years later, in \cite{Cha} M. Chaperon said that no infinite-dimensional version of Holmgren's theorem seems to be known.

The main purpose of this paper is to establish Cauchy-Kowalevski and Holmgren type theorems on $\mathbb{R}^{\infty}$, which is a countable Cartesian product of $\mathbb{R}$.
In some sense, this is quite natural because $\mathbb{R}^{\infty}$ is an infinite dimensional counterpart of $\mathbb{R}^{n}$, the $n$-fold  Cartesian product of $\mathbb{R}$.
Since there are various different topologies on $\mathbb{R}^{\infty}$, we have more freedom and flexibility to introduce suitable assumptions on equations under consideration. Nevertheless, a large family of local derivatives of functions on $\mathbb{R}^{\infty}$ make the analysis in our work more complicated.
Clearly, we cannot apply the known abstract Cauchy-Kowalevski and Holmgren type theorems (in the literatures) in the setting of Banach spaces to obtain our results. Indeed, our working space, $\mathbb{R}^{\infty}$, is NOT a Banach space!

There exists many different definitions of analyticity in infinite dimensions. As far as we know, the concept of analyticity for functions of infinitely many variables began with H.~von Koch in 1899 (\cite{Koc}). H.~von Koch (\cite{Koc}) introduced a monomial approach to holomorphic functions on infinite dimensional polydiscs which was further developed by D.~Hilbert in 1909 (\cite{Hil}). After the pioneering work of H.~von Koch and D.~Hilbert, it is clear from M.~Fr\'echet \cite{Fre1, Fre2} and R.~G\^{a}teaux \cite{Gat1, Gat2} that the power series expansion in terms of homogeneous polynomials seems more suitable for  the analyticity in the context of Banach space (e.g., \cite{Din} for the extensive works on infinite dimensional complex analysis starting from 1960s). Nevertheless, in the recent decades,  it was found that Hilbert's definition of analyticity is also useful in some problems. For example, in 1987, R.~Ryan (\cite{Rya}) discovered that every entire function $f$ on $\ell^1$ (the usual Banach space of absolutely summable sequences of real numbers) has a monomial expansion which converges and coincides with $f$; in 1999, L.~Lempert (\cite{Lem}) proved this holds for any open ball of $\ell^1$; while in 2009, A.~Defant, M.~Maestre and C.~Prengel (\cite{DMP}) showed that the monomial expansion of any holomorphic function on the Reinhardt domain $\cal R$ in a Banach sequence space converges uniformly and absolutely on any compact subsets of $\cal R$. In this work, we shall use the definition of analyticity introduced by H.~von Koch and D.~Hilbert.

We shall modify the classical method of majorants to derive a Cauchy-Kowalevski type theorem in $\mathbb{R}^{\infty}$. Based on this result, we then employ some tools from abstract Wiener spaces (developed by L.~Gross in \cite{Gro1}) and especially an infinite dimensional divergence theorem (\cite{Goo}) to establish a Holmgren type theorem with infinite number of variables. For the later, the basic idea is more or less the same as that in finite dimensions but the technique details are much more complicated. Indeed, it is well-known that, compared with its finite dimensional counterpart, the analysis tools in infinite dimensions are much less developed.

The rest of this paper is organized as follows. Section 2 is of preliminary nature, in which we first introduce suitable topologies for $\mathbb{R}^{\infty}$ and for $\mathbb{R}^{\infty}$ plus an infinity point; then, we give the definition of analyticity for functions of infinitely many variables; also, we present a brief introduction to the theory of abstract Wiener space and a divergence theorem. In Section 3, we show a Cauchy-Kowalevski type theorem of infinitely many variables. Section 4 is devoted to establishing our Holmgren type theorem.

We refer to \cite{YuZh1} for the details of proofs of the results announced in this paper and some other results in
this context.

\section{Preliminaries }

\subsection{A family of topologies}

There exists only one useful topology on finite dimensional space $\mathbb{R}^{n}$. However, we need to use a family of topologies on  $\mathbb{R}^{\infty}$. 

Denote by $\mathscr{B}_{\infty}$ the class of sets (in $\mathbb{R}^{\infty}$):
$(x_i)+B_r^{\infty}\triangleq \{(x_i+y_i):\, (y_i)\in B_r^{\infty}\}$,
 where $(x_i)\in \mathbb{R}^{\infty}$, $r\in(0,+\infty)$ and $B_r^{\infty}\triangleq\{(y_i)\in \mathbb{R}^{\infty}:\,\sup_{1\leq i< \infty}|y_i|<r\}.$ Then $\mathscr{B}_{\infty}$ is a base for a topological space, denoted by $\mathscr{T}_{\infty}$.
For any $(x_i),(y_i)\in \mathbb{R}^{\infty}$, define $d_{\infty}((x_i),(y_i))\triangleq \min \{1,\sup_{1\leq i<\infty}|x_i-y_i|\}$. Then, it is easy to show the following result:

\begin{proposition}\label{prop 3}
The following assertions hold:
\begin{itemize}
\item[(1)]  $(\mathbb{R}^{\infty},\mathscr{T}_{\infty})$ is not a topological vector space;
\item[(2)]  $(\mathbb{R}^{\infty},\mathscr{T}_{\infty})$ is compatible with the metric $d_{\infty}$ and $(\mathbb{R}^{\infty},d_{\infty})$ is a nonseparable complete metric space.
\end{itemize}
\end{proposition}

For any $p\in[1,\infty)$ and $(x_i),(y_i)\in \mathbb{R}^{\infty}$, define
$\text{d}_p((x_i),(y_i))\triangleq \min \{1, (\sum_{i=1}^{\infty}|x_i-y_i|^p )^{\frac{1}{p}} \}.$ Then, $(\mathbb{R}^{\infty}, d_p)$ is a complete metric space. Denote by $\mathscr{T}_p$ the topology induced by the metric $d_p$, and by $\mathscr{B}_p$ the class of sets: $(x_i)+B_r^{p}$, where $(x_i)\in \mathbb{R}^{\infty}$, $r\in(0,+\infty)$ and $B_r^{p}\triangleq \{(x_i)\in \mathbb{R}^{\infty}:\, (\sum_{i=1}^{\infty}|x_i|^p )^{\frac{1}{p}}<r \}$. Obviously, $\mathscr{B}_p$ is a base for the topology space $\mathscr{T}_p$. Also, denote by $\ell^p$ ({\it resp.} $\ell^\infty$) the usual Banach space of sequences $(x_i)\in \mathbb{R}^{\infty}$ so that $(\sum_{i=1}^{\infty}|x_i|^p )^{\frac{1}{p}}<\infty$ ({\it resp.} $\sup_{1\leq i< \infty}|x_i|<\infty$).

Likewise, for any $p\in(0,1)$ and $(x_i),(y_i)\in \mathbb{R}^{\infty}$, define $\text{d}_p((x_i),(y_i))\triangleq \min \{1,  (\sum_{i=1}^{\infty}|x_i-y_i|^p ) \}$.
Then $d_p$ is a metric on $\mathbb{R}^{\infty}$. Denote by $\mathscr{T}_p$ the topology induced by the metric $d_p$, and by  $\mathscr{B}_p$ the class of sets: $(x_i)+B_r^{p}$, where $(x_i)\in \mathbb{R}^{\infty}$, $r\in(0,+\infty)$ and $B_r^{p}\triangleq \{(x_i)\in \mathbb{R}^{\infty}:\,\sum_{i=1}^{\infty}|x_i|^p<r \}$.
Similarly to Proposition \ref{prop 3}, $(\mathbb{R}^{\infty},\mathscr{T}_{p})$ and $\mathscr{B}_p$ is a base for the topology space $\mathscr{T}_p$ when $0< p< \infty$.

     Denote by $\mathscr{T}$ the (usual) product topology on $\mathbb{R}^{\infty}$. Now we have defined a family of topologies on $\mathbb{R}^{\infty}$. The following result  shows some relations between these topologies.

\begin{proposition} \label{toplogies inclusion proposition}
The inclusion relations $\mathscr{T}\subsetneqq \mathscr{T}_{\infty}\subsetneqq \mathscr{T}_q\subsetneqq\mathscr{T}_p$ hold for any $0<p<q<\infty$.
\end{proposition}

Let $\widetilde{\mathbb{R}^{\infty} }\triangleq \mathbb{R}^{\infty} \sqcup\{\infty\}$, where $\infty$ is any fixed point not belonging in $\mathbb{R}^{\infty}$. We consider the following family of sets (in $\widetilde{\mathbb{R}^{\infty} }$):
$\mathscr{B}^p\triangleq \mathscr{B}_p \bigcup \Big\{\big\{(x_n)\in  \mathbb{R}^{\infty} :x_n\neq 0 \text{ for each }n\in \mathbb{N}\,\text{ and } \,\,\sum_{n=1}^{\infty}\frac{1}{|x_n|^p}<r \big\}\sqcup\{\infty\}:\,r>0\Big\}$ for $0<p<\infty$, and $\mathscr{B}^{\infty}\triangleq \mathscr{B}_{\infty} \bigcup \Big\{\big\{(x_n)\in \mathbb{R}^{\infty} :x_n\neq 0 \text{ for each }n\in \mathbb{N}\,\text{ and } \,\,\sup_{1\leq n<\infty}\frac{1}{|x_n|}<r \big\}\sqcup \{\infty\}:\,r>0\Big\}.$ Then $\mathscr{B}^p$ is a base for a topological space on $\widetilde{\mathbb{R}^{\infty} }$ which we will denoted by $\mathcal {T}^p$ for $0<p\leq \infty$. It is easily seen that the subspace topology of $(\widetilde{\mathbb{R}^{\infty} },\mathcal {T}^p)$ on $\mathbb{R}^{\infty}$ is $\mathscr{T}_p$ but $\mathcal {T}^p$  is not the usual one-point compactification of $\mathscr{T}_p$ (and therefore we use the notion $\widetilde{\mathbb{R}^{\infty} }$ rather than $\widehat {\mathbb{R}^{\infty} }$ for the usual one-point compactification).

\subsection{Analyticity for Functions of Infinitely Many Variables}
We denote by $\mathbb{N}^{(\mathbb{N})}$ the set of all finitely supported sequences of nonnegative integers. For $\textbf{x}=(x_i)\in \mathbb{R}^{\infty}$, $\alpha=(\alpha_i)\in \mathbb{N}^{(\mathbb{N})}$ with $\alpha_k=0$ for $k\geq n+1$ and some $n\in \mathbb{N}$, write $\textbf{x}^{\alpha}\triangleq x_1^{\alpha_1}x_2^{\alpha_2}\cdots x_n^{\alpha_n}$, which is called a monomial on $\mathbb{R}^{\infty}$. The following definition of analyticity is essentially from \cite{Koc} and \cite{Hil}.
\begin{definition}
Suppose $f$ is a real-valued function defined on a subset $D$  of $\mathbb{R}^{\infty}$ and $\textbf{x}_0\in D$. If for each $\alpha\in \mathbb{N}^{(\mathbb{N})}$, there exist a real number $c_{\alpha}$ (depending on $f$ and $x_0$) such that the series $\sum_{\alpha \in \mathbb{N}^{(\mathbb{N})}}  c_{\alpha}\textbf{x}^{\alpha}$ is absolutely convergent for some $\textbf{x}=(x_i)\in \mathbb{R}^{\infty}$ with $x_i\neq 0$ for each $i\in \mathbb{N}$, and
for each $\textbf{h}=(h_i)\in \mathbb{R}^{\infty}$ with $|h_i|\leq |x_i|$ for all $i\in\mathbb{N}$, it holds that $\textbf{x}_0+\textbf{h}\in D$ and
\begin{eqnarray}\label{monomial expansion}
f(\textbf{x}_0+\textbf{h})=\sum_{\alpha \in \mathbb{N}^{(\mathbb{N})}}  c_{\alpha}\textbf{h}^{\alpha},
\end{eqnarray}
then  $f$ is called analytic near $\textbf{x}_0$ (with the monomial expansion (\ref{monomial expansion})). In this case, we write $D_f^{\textbf{x}_0}\triangleq \{\textbf{h}\in \mathbb{R}^{\infty}:\sum_{\alpha \in \mathbb{N}^{(\mathbb{N})}}  |c_{\alpha}\textbf{h}^{\alpha}|<\infty \}$, and call the set $D_f^{\textbf{x}_0}$ the convergence domain of monomial expansion (\ref{monomial expansion}).
\end{definition}
\begin{definition}
Suppose $f$ is a real-valued function defined on a subset $D$ of $\mathbb{R}^{\infty}$, $\textbf{x}_0\in D$ and $0<p<\infty$. If there exists $\textbf{x}=(x_i)\in\mathbb{R}^{\infty}$ with $\sum_{i=1}^{\infty}\frac{1}{|x_i|^p}<\infty$ such that $f$ is analytic near $\textbf{x}_0$ with the monomial expansion (\ref{monomial expansion}) and $\sum_{\alpha \in \mathbb{N}^{(\mathbb{N})}}  |c_{\alpha}\textbf{x}^{\alpha}|<\infty$, then the monomial expansion of $f$ near $\textbf{x}_0$ is called absolutely convergent at a point near $\infty$ in the topology $\mathcal {T}^{p}$.
\end{definition}
\begin{example}
(\textbf{Riemann Zeta Function}) Recall that the Riemann zeta function is defined by $\zeta(s)=\sum_{n=1}^{\infty}\frac{1}{n^s}$ when $s\in \mathbb{C}$ with Re$(s)>2$. Suppose that $\{p_n\}_{n=1}^{\infty}$ is the collection of all positive prime numbers. Then for $s\in \mathbb{R}$ with $s>2$ we have $\zeta(s)=\prod_{n=1}^{\infty}\frac{1}{1-p_n^{-s}}
=\sum_{\alpha=(\alpha_1,\cdots,\alpha_n) \in \mathbb{N}^{(\mathbb{N})}}p_1^{-s\alpha_1}\cdots p_n^{-s\alpha_n}
=\sum_{\alpha \in \mathbb{N}^{(\mathbb{N})}}(p_1^{-s},\cdots,p_n^{-s},\cdots)^{\alpha}
=f(p_1^{-s},\cdots,p_n^{-s},\cdots),$ where $f(\textbf{x})\triangleq\sum_{\alpha \in \mathbb{N}^{(\mathbb{N})}}\textbf{x}^{\alpha}$ with definition domain $D=\{\textbf{x}=(x_i)\in\mathbb{R}^{\infty}:\,\sum_{i=1}^{\infty}|x_i|<1\}$ is the function of geometric series of infinitely many variables.
\end{example}


Naturally, the notion of majority function is as follows:

\begin{definition}
Suppose $f$ and $F$ are analytic near $\textbf{x}_0\in \mathbb{R}^{\infty}$ with monomial expansions $f(\textbf{x})=\sum_{\alpha \in \mathbb{N}^{(\mathbb{N})}} c_{\alpha}(\textbf{x}-\textbf{x}_0)^{\alpha}$ and $F(\textbf{x})=\sum_{\alpha \in \mathbb{N}^{(\mathbb{N})}} C_{\alpha}(\textbf{x}-\textbf{x}_0)^{\alpha}$, respectively,
where $c_{\alpha}\in\mathbb{R}$ and $C_{\alpha}\geq 0$ for each $\alpha \in \mathbb{N}^{(\mathbb{N})}$. If $|c_{\alpha}|\leq C_{\alpha}$ for all $\alpha \in \mathbb{N}^{(\mathbb{N})}$, then $F$ is called a majority function of $f$ near $\textbf{x}_0$.
\end{definition}

Suppose $B$ is a real Banach space and $U$ is an open set of $B$. For each $n\in \mathbb{N}$, a function $g$ from $B$ into $\mathbb{R}$ is called continuous $n$-homogeneous polynomial if there exists a continuous $n$-linear map $T$ from $\prod_{i=1}^{n}B$ into $\mathbb{R}$ such that $g(x)=T(x,\cdots,x)$ for each $x\in B$. For $n=0$, we call any function from $B$ into $\mathbb{R}$ with constant value a continuous $0$-homogeneous polynomial.
A function $f$ (defined on a subset $D$  of $B$) is called analytic on some subset $U\subset D$ if for each $\xi\in U$ there exist a sequence $\{P_nf(\xi)\}_{n=0}^{\infty}$ of continuous $n$-homogeneous polynomials  on $B$ and a radius $r>0$ such that $\xi+B(r)\subset U$ and $f(\xi+x)=\sum_{n=0}^{\infty} P_nf(\xi)(x)$ uniformly in $B(r)$, where $B(r)\triangleq \{y\in B: ||y||<r\}$. It is easy to see that, the function $g(\cdot)\triangleq f(\xi+\cdot)$ is Fr\'{e}chet differentiable in $B(r)$ (See \cite{Din} for more details).

We emphasize that $\mathbb{R}^{\infty}$ is not a norm space. Hence, the analyticities by monomial expansions and by power series expansions are two distinct notions.
Nevertheless, for every $\textbf{x}\in\mathbb{R}^{\infty}$, $r\in(0,1)$ and $p\in[1,+\infty]$, it holds that $\textbf{x}+B_r^{p}\subset \mathbb{R}^{\infty}$. Motivated by this simple observation, we have the following definition.
\begin{definition}
Suppose $p\in[1,\infty],\,U\subset \mathbb{R}^{\infty}$, $f$ is a real-valued function defined on $U$ and $\textbf{x}\in U$. We call $f$ is Fr\'{e}chet differentiable with respect to $\ell^p$ in a neighborhood of $\textbf{x}$ in the topology $\mathscr{T}_p$ if there exists $r\in(0,1)$ such that $\textbf{x}+B_r^{p}\subset U$, and the function defined by $g(\cdot)\triangleq f(\textbf{x}+\cdot)$ is Fr\'{e}chet differentiable with respect to the Banach space $\ell^p$ in $B_r^{p}$.
\end{definition}

\subsection{Abstract Wiener Space and Derivatives}

Let us recall the notion of abstract Wiener space, which will be of crucial importance in the proof of the Holmgren type theorem. The materials in this subsection are from \cite{Dri, Fer, Gro, Gro1}.

Suppose $X$ is a real separable Banach space, and denote by $X^*$ its dual space and by $\langle \cdot,\cdot\rangle$ the natural pairing from $X^*\times X$ into $\mathbb{R}$. Let $Y$ be another
Banach space. Denote by ${\cal L}(X; Y)$ the
Banach space of all bounded linear operators
from $X$ to $Y$, with the usual operator
norm (When $Y=X$, we simply write ${\cal L}(X)$
instead of ${\cal L}(X; Y)$). A subset $C$ of $X$ is called a cylinder set if it is of the form:
$$C=\big\{x\in X:(\langle y_1,x\rangle,\cdots,\langle y_n,x\rangle)\in E\big\},$$
where $n\in\mathbb{N}$, $y_1,\cdots,y_n\in X^*$ and $E$ is a Borel set in $\mathbb{R}^n$. If $L$ is a finite-dimensional subspace of $X^*$ such that $\{y_1,\cdots,y_n\}\subset L$, then $C$ is said to be based on $L$. Clearly, the collection of cylinder sets in $X$ is an algebra $\mathscr{R}$ and the collection of cylinder sets based on $L$ is a $\sigma$-algebra which will be denoted by $\mathscr{S}_L$. We call a nonnegative set function $\mu$ on $\mathscr{R}$ a cylinder set measure on $X$ if $\mu(X)=1$ and $\mu$ is countably additive on $\mathscr{S}_L$ for each finite-dimensional subspace $L$ of $X^*$.

Suppose $H$ is a real separable Hilbert space with inner product $(\cdot,
\cdot)$ and norm $|\cdot|=\sqrt{(\cdot,\cdot)}$. Then every cylinder set of $H$ is of the form
$$C=\big\{x\in H:Px\in E\big\},$$ where $P$ is a finite-dimensional projection in $H$ and $E$ is a Borel set in $PH$. For any $t>0$, a (typical) cylinder set measure $\mu_t$ is defined by
$$\mu_t(C)\triangleq\frac{1}{(2\pi t)^{-\frac{n}{2}}}\int_E e^{-\frac{|x|^2}{2t}}\,\mathrm{d}x,$$
where $C$, $P$ and $E$ are given above, $n=\dim PH$, and $\mathrm{d}x$ is the Lebesgue measure in $PH$.
A measurable semi-norm on $H$ is a semi-norm $||\cdot||$ on $H$ with the property that for every $\epsilon>0$, there exists a finite-dimensional projection $P_0$ such that for every finite-dimensional projection $P$ orthogonal to $P_0$ it holds that $\mu_1(\{x\in H: ||Px||>\epsilon\})<\epsilon.$ As a consequence of the definition of measurable semi-norm, every measurable semi-norm $||\cdot||$ is dominated by the Hilbert norm, i.e., there exists a constant $C$ such that $||x||\leq C|x|$ for all $x\in H$. If $||\cdot||$ is a measurable norm, we denote by $B$ the completion of $H$ with respect to $||\cdot||$. Then $B$ is a separable  Banach space. There is a natural embedding $i$ from $H$ into $B$ whose image is dense in $B$. Then $i^*$ is also an embedding from $B^*$ into $H^*$. Since $H^*$ can be identified with $H$ we have the following inclusion relations $B^*\subset H\subset B.$ Furthermore, we should note that an element $x$ of $H$ is in $B^*$ if and only if there exists a constant $C>0$ such that the inequality $|(x,y)|\leq C ||y||$ holds for all $y$ in $H$. Then $\mu_t$ induces a cylinder set measure $m_t$ in $B$ as follows. If $y_1,\cdots,y_n\in B^*$ and $E$ is a Borel set of $\mathbb{R}^n$, we define
$$m_t(\{x\in B:(\langle y_1,x\rangle,\cdots,\langle y_n,x\rangle)\in E\})\triangleq \mu_t(\{x\in H:((y_1,x),\cdots,( y_n,x))\in E\}).$$
In \cite{Gro1}, it was proved that $m_t$ is countably additive on the cylinder sets of $B$. By the Carath\'{e}odory extension theorem, it can be uniquely extended to the 
Borel sets of $B$ as a measure, denoted by $p_t$. The triple $(i,H,B)$ is called an abstract Wiener space and the measure $p_t$ is called the Wiener measure on $B$ with variance parameter $t$. For any $x\in B$ and Borel subset $A$ of $B$, we define $p_t(x,A)\triangleq p_t(A-x).$ By \cite{Gro}, one has the following formula for the Wiener measure.
\begin{proposition}
For any $s,t\in(0,+\infty)$ and $x,y\in B$, $p_t(x,\cdot)$ and $p_s(y,\cdot)$ are equivalent measures if and only if $s=t$ and $x-y\in H$. Otherwise they are mutually singular. Furthermore, it holds that $p_{ts}(A)=p_t(s^{-\frac{1}{2}}A)$ for any Borel subset $A$ of $B$.
\end{proposition}
In terms of the inclusion relations $B^*\subset H\subset B$, one has a useful product decomposition of Wiener measure $p_t$ as follows. Suppose $K$ is a finite-dimensional subspace of $B^*$ and $L$ is its annihilator in $B$. If $\{y_1,\cdots,y_n\}$ is an orthonormal basis of $K$ then we define $Qx=\sum_{j=1}^{n}\langle y_j,x\rangle y_j,\,x\in B.$ Then $Q$ is a continuous linear operator from $B$ into $B$. Obviously, the range of $Q$ is $K$ and null space of $Q$ is $L$. We thus get $B=K\oplus L$. Let $K^{\perp}$ be the orthogonal complement of $K$ in $H$. Then $K^{\perp}\subset L$ and $L$ is the closure in $B$ of $K^{\perp}$. It is easy to check that the restriction of a measurable norm to a closed subspace is again a measurable norm. Thus there is a Wiener measure $p_t^{'}$ on the space $L$. Let $\mu_t^{'}$ denote the typical Gaussian measure in $K$ then in the Cartesian product decomposition $B=K\times L$ there holds $p_t=\mu_t^{'}\times p_t^{'}.$

The following exponential integrability of Wiener measure, discovered in \cite{Fer}, is very useful.
\begin{theorem}\label{Fernique's Theorem}
{\rm (\textbf{Fernique's Theorem})} For any fixed $t\in(0,\infty)$, there exists $\epsilon=\epsilon(p_t)>0$ such that
\begin{eqnarray*}
\int_B e^{\epsilon ||x||^2}\,p_t(\mathrm{d}x)<\infty.
\end{eqnarray*}
\end{theorem}

Denote by $\mathcal{B}(B)$ the collection of Borel sets of $B$. The following density result (\cite[Proposition 39.8]{Dri}) will also be used later.
\begin{proposition}\label{prop 39.8 in Drive}
Suppose $\mu$ is a probability measure on $\mathcal{B}(B)$ so that for every $\varphi\in B^*$, there exists a constant $\epsilon=\epsilon( ||\varphi||_{B^*})>0$ such that $\int_Be^{\epsilon |\varphi|}\,\mathrm{d}\mu<\infty$ (Here, $|\varphi|$ stands for the absolute value of $\varphi$). Then $\mathcal {F}\triangleq\{P(\varphi_1,\cdots,\varphi_n):n\in\mathbb{N},\varphi_i\in B^*,1\leq i\leq n, P\text{ is a real polynomial with $n$ variables }\}$ is dense in $L^p(B,\mathcal {B}(B),\mu)$ for $1\leq p<\infty$.
\end{proposition}

In what follows, $(i,H,B)$ is a given abstract Wiener space. For a real valued function $f$ defined on an open set $O$ of $B$, one has two definitions of derivatives. Firstly, for $x\in O$ if $f$ is $B$-Fr\'{e}chet differentiable at $x$, we will denote the corresponding derivative at $x$ by $f'(x)$. Secondly, we have a function defined by $g(h)=f(x+ih)$ which is a function on a neighborhood of $0$ in $H$. If $g$ is $H$-Fr\'{e}chet differentiable at $0$, we will denote the corresponding derivative at $0$ by $Df(x)$. From \cite{Gro} we see that the second derivative is weaker than the first derivative. The second derivative is sometimes called Malliavin derivative, which will be used in the proof of the Holmgren type theorem in this paper. For $n\geq 2$, we will use the notations $f^{n}(x)$ and $D^{n}f(x)$ to denote the corresponding higher order derivatives.
\subsection{A Divergence Theorem in Abstract Wiener Space}
In this subsection, we will give a brief exposition of surface measures and a divergence theorem  in the abstract Wiener space $(i,H,B)$, developed in \cite{Goo}.

Firstly, we introduce the concepts of ``smooth" functions and surfaces.
\begin{definition}
A real-valued function $g$ defined on an open subset $U$ of $B$ is called an $H$-$C^1$ function if $g$ is continuous and $H$-Fr\'{e}chet  differentiable on $U$, the map $Dg:U\to H^*$ is continuous and the vector $Dg(x)\in B^*$ for each $x\in U$.
\end{definition}
\begin{definition}
A subset $S$ ({\it resp.} $V$) of $B$ is called an $H$-$C^1$ surface ({\it resp.} to have an $H$-$C^1$ boundary $\partial V$) if for each $x\in S$ ({\it resp.} $x\in\partial V$) there is an open neighborhood $U$ of $x$ in $B$ and an $H$-$C^1$ function $g$ defined on $U$ such that $Dg(x)\neq 0$ and $S\cap U=\{y\in U: g(y)=0\}$ ({\it resp.} $V\cap U=\{y\in U:g(y)<0\}$).
\end{definition}
Secondly, we introduce the notion of local coordinate for the above $H$-$C^1$ surface $S$. We begin with the concept of normal projection:
\begin{definition}
A one dimensional orthogonal projection $P$ on $H$ is called a normal projection for $S$ at $x\in S$ if $|P(y-x)|=o(|y-x|)$ as $|y-x|\to 0$ for all $y\in S$ so that $y-x\in H$.
\end{definition}
Denote by $I$ the identity operator on $H$. One can show the following result:
\begin{proposition}\label{prop 2 in Godman DT}
There exists a unique map $N_{\cdot}:S\to {\cal L}(B;B^*)$ such that
\begin{itemize}
 \item[(1)] For each $x\in S$ the restriction of $N_x$ to $H$ is a normal projection for $(S-x)\cap H$ at $0$;
\item[(2)] For each $x$ the map $J_x=I-N_x$ is a homeomorphism of an open neighborhood of $x$ in $S$ onto an open subset in the null space of $N_x$;
\item[(3)] The map $N:S\to {\cal L}(H)$ is continuous.
\end{itemize}
\end{proposition}

For any $y\in S$, by Proposition \ref{prop 2 in Godman DT}, there is an element $h$ in $B^*$ with $|h|=1$ and a neighborhood $W$ of $y$ in $S$ such that
$N_yh=h$, $|N_wh|>0$ for all $w$ in $W$, and $J_y=I-N_y$ is a homeomorphism of $W$ into the null space of $N_y$.
We call the above element $h$ {\it a unit normal vector} at $y$, and $W$ {\it a coordinate neighborhood} of $y$.
For a fixed $x$ in $B$ and $t>0$ we define a measure $\rho(t,x,W,\cdot)$ on the Borel sets of $W$ by
$$\rho(t,x,W,E)\triangleq\frac{1}{\sqrt{2\pi t}}\int_{J_y(E)} \frac{1}{|N_{J_y^{-1}z}h|}\exp\Big[-\frac{|N_y(J_y^{-1}z-x)|^2}{2t}\Big]\,p_t^{'}(J_y x,\mathrm{d}z),$$
where $p_t^{'}$ is the Wiener measure on the null space of $N_y$. We call $\rho(t,x,W,\cdot)$ {\it a local version of normal surface measure} with dilation parameter $t$ and translation variable $x$.

One has the following result:

\begin{theorem}\label{theorem 1 in Goodman DT}
 For any $x\in B$ and $t>0$, there is a unique measure $\sigma_t(x,\cdot)$ on the Borel sets of $S$ such that for any local version of normal surface measure  $\rho(t,x,W,\cdot)$ on a coordinate neighborhood $W$ in $S$ and any Borel subset $E$ of $W$, it holds that $\sigma_t(x,E)=\rho(t,x,W,E)$.
\end{theorem}

The measure $\sigma_t(x,\cdot)$ given in Theorem \ref{theorem 1 in Goodman DT} is called {\it a normal surface measure} on $S$ with dilation parameter $t$ and translation variable $x$.

In the sequel, $V$ is a given nonempty open set of $B$ and has an $H$-$C^1$ boundary $\partial V$. We also need the following two concepts.
\begin{definition}
A map ${\bf n}:\partial V\to H^*$  is called {\it an outward normal map} for $V$ provided that ${\bf n}(y)$ is a unit normal vector at $y$ for the surface $\partial V$ and $y-s{\bf n}(y)\in V$ for any small $s>0$.
\end{definition}

\begin{definition}
An ($H$-valued) function $F$ defined on an open subset $U$ of $B$ is called to have finite divergence at $x\in U$ if $F$ is $H$-Fr\'{e}chet differentiable at $x$ and $DF(x)$ is an operator of trace class on $H$. For such a function $F$, the divergence of $F$ at $x$ is defined by the  trace of $DF(x)$, and denoted by $(\text{div }F)(x)$.
\end{definition}
The following result will play a key role in the proof of our Holmgren type theorem of infinite many variables.
\begin{theorem}\label{Divergence Theorem}{\rm (\textbf{Divergence Theorem})} Assume that $F:V\cup \partial V\to H$ is a continuous function with finite divergence on $V$ and that $F$ is uniformly bounded with respect to the $B^*$-norm on $V$. If for some $x\in B$ and $t>0$, the function $|F(\cdot)|$ is $\sigma_t(x,\cdot)$-integrable on $\partial V$ and the trace class operator norm of $DF$ is $p_t(x,\cdot )$-integrable on $V$, then
\begin{eqnarray*}
\int_V \bigg[({\rm div }\;F)(y)-\frac{\langle F(y),y-x\rangle}{t}\bigg]\,p_t(x,\mathrm{d}y)=\int_{\partial V}\langle F(y), {\bf n}(y)\rangle\,\sigma_t(x,\mathrm{d}y).
\end{eqnarray*}
\end{theorem}

\section{Cauchy-Kowalevski Type Theorem  of Infinitely Many Variables}
This section is devoted to a study of the following form of initial value problem:
\begin{equation}\label{one order quasi-linear partial differential equation}
\left\{
\begin{array}{ll}
 \displaystyle\partial^m_t u(t,\textbf{x})=f\big(t,\textbf{x},u, \partial^{\beta}_{\textbf{x}}\partial^{j}_t  u \big),\\[2mm]\displaystyle
 \partial^k_tu(t,\textbf{x})\mid_{t=0}=\phi_k(\textbf{x}),\ \ k=0,1,\cdots,m-1.
 \end{array}
\right.
\end{equation}
Here $m$ is a given positive integer, $t\in \mathbb{R}, \textbf{x}=(x_i)\in \mathbb{R}^{\infty}$, $\partial^{\beta}_{\textbf{x}}\;{\buildrel \triangle \over =}\;\partial^{\beta_1}_{x_1^{\beta_1} }\cdots\partial^{\beta_n}_{x_n^{\beta_n} }$  for $\beta=(\beta_1,\cdots,\beta_n,0,\cdots)\in \mathbb{N}^{(\mathbb{N})}$ (for some positive $n\in\mathbb{N}$), the unknown $u$ is a real-valued function depending on $t$ and $\textbf{x}$; $f$ is a non-linear real-valued function depending on $t,$ $\textbf{x}$, $u$ and all of its derivatives of the form $\partial^{\beta}_x\partial^{j}_t  u,\,\beta\in \mathbb{N}^{(\mathbb{N})}, j<m, 1\leq|\beta|+j\leq m.$ Note that the values $u(0,(0)), \partial^{\beta}_x\partial^{j}_t  u(0,(0)), \beta\in \mathbb{N}^{(\mathbb{N})},\,j<m,\,1\leq|\beta|+j\leq m$ are determined by (\ref{one order quasi-linear partial differential equation}). For simplicity, we write these determined values by $u(0,(0))=u_0,\partial^{\beta}_x\partial^{j}_t  u(0,(0))=w_{\beta,j}^0,\beta\in \mathbb{N}^{(\mathbb{N})},\,j<m,\,1\leq|\beta|+j\leq m$ and $\textbf{w}^0= (w_{\beta,j}^0)_{\beta\in \mathbb{N}^{(\mathbb{N})},\,j<m,\,1\leq|\beta|+j\leq m}.$ One can see that $f$ is a function on $\mathbb{R}\times \mathbb{R}^{\infty}\times \mathbb{R}\times \mathbb{R}^{\infty}$ which is also a countable Cartesian product of $\mathbb{R}$. Therefore, we may identify $\mathbb{R}\times \mathbb{R}^{\infty}\times \mathbb{R}\times \mathbb{R}^{\infty}$ with $\mathbb{R}^{\infty}$.  We suppose that $f$ is analytic near $(0,(0),u(0,(0)), \partial^{\beta}_x\partial^{j}_t  u(0,(0)) )$ with the monomial expansion $f\big(t,(x_i),u, \partial^{\beta}_x\partial^{j}_t  u \big)=\sum_{\alpha\in\mathbb{N}^{(\mathbb{N})}}C_{\alpha}\big(t,(x_i),u-u_0, \partial^{\beta}_x\partial^{j}_t  u-w_{\beta,j}^0 \big)^{\alpha}$
and let $F(t,\textbf{x},u, \textbf{w})\triangleq\sum_{\alpha\in \mathbb{N}^{(\mathbb{N})}}|C_{\alpha}| (t,\textbf{x},u, \textbf{w})^{\alpha},$ where
$\textbf{w}= (w_{\beta,j})_{\beta\in \mathbb{N}^{(\mathbb{N})},\,1\leq |\beta|+j\leq m,\,j<m}\in \mathbb{R}^{\infty}$ and $\textbf{x}= (x_i)\in \mathbb{R}^{\infty}.$ We also suppose that for each $0\leq k\leq m-1$, $\phi_k$ is analytic near $(0,(0))$ with monomial expansion $\phi_k(\textbf{x})=\sum_{\alpha\in \mathbb{N}^{(\mathbb{N})}}C_{\alpha,k}\textbf{x}^{\alpha}$ and let $\Phi_k (t,\textbf{x},u, \textbf{w})\triangleq\sum_{\alpha\in \mathbb{N}^{(\mathbb{N})}}|C_{\alpha,k}|\textbf{x}^{\alpha}.$

By means of the majorant method, we can show the following Cauchy-Kowalevski type theorem of infinitely many variables:
\begin{theorem}\label{Infinite Dimensional Cauchy-Kowalevski Theorem}
Suppose $1\leq p<\infty$ and that the monomial expansions of $\Phi_k$, $0\leq k\leq m-1$ and $F$ near $(0,(0))$ are absolutely convergent at a point near $\infty$ in the topology $\mathcal {T}^p$. Then the Cauchy problem (\ref{one order quasi-linear partial differential equation}) admits a locally analytic solution (near $(0)$), which is unique in the class of analytic functions under the topology $\mathscr{T}_{p'}$ where $p'$ is the usual H\"older conjugate of $p$. Furthermore, the solution $u$ is Fr\'{e}chet differentiable with respect to $\ell^{p'}$ in a neighborhood of $(0)$ in the topology $\mathcal {T}^{p'}$ and the corresponding Fr\'{e}chet derivative  $Du$ is continuous.
\end{theorem}

\begin{example}
(\textbf{The Schr\"{o}dinger operator of infinitely many number of particles}) This example is from \cite{Bere86}. Suppose that $\{a_k\}_{k=1}^{\infty}$ is a sequence of nonnegative real numbers. The following operator
\begin{eqnarray}\label{operator in QFT}
(Lu)(\textbf{x})\triangleq-\frac{1}{2}\sum_{k=1}^{\infty}a_k\bigg(\frac{\partial^2 u(\textbf{x})}{\partial x_k^2}-2x_k\frac{\partial u(\textbf{x})}{\partial x_k}\bigg),\ \ \textbf{x}=(x_k)\in \mathbb{R}^{\infty},
\end{eqnarray}
for $u\in \mathcal {P}(\mathbb{R}^{\infty})$ which is the space of cylindrical polynomials on $\mathbb{R}^{\infty}$, i.e., polynomials depending only on finitely many variables. Note that in quantum mechanics the operator $N_k\triangleq-\frac{1}{2}a_k\Big(\frac{\partial^2  }{\partial x_k^2}-2x_k\frac{\partial }{\partial x_k}\Big)$ is an operator of energy of one-dimensional harmonic oscillator with unit mass and frequency in the space of $L^2\Big(\mathbb{R},\frac{e^{-x_k^2}}{\sqrt{\pi}}\,\mathrm{d}x_k\Big)$. Thus the operator (\ref{operator in QFT}) describes a system consisting infinitely many noninteracting oscillators with frequency $a_k\geq 0,k\in \mathbb{N}$.  If $a_1=1$ and we let $t=x_1$ then the equation $Lu=0$ is equivalent to the following one:
$$\frac{\partial^2 u(t,\textbf{x})}{\partial t^2}=2t\frac{\partial u(t,\textbf{x})}{\partial t}-\sum_{k=2}^{\infty}a_k\Big(\frac{\partial^2 u(t,\textbf{x})}{\partial x_k^2}-2x_k\frac{\partial u(t,\textbf{x})}{\partial x_k}\Big),
$$
which is of form (\ref{one order quasi-linear partial differential equation}). One can easily see that the above equation satisfies the assumption in Theorem \ref{Infinite Dimensional Cauchy-Kowalevski Theorem} if and only if there exist $(b_k),(c_k)\in \mathbb{R}^{\infty}$ such that $\sum_{k=2}^{\infty}\big(\frac{1}{|b_k|^p}+\frac{1}{|c_k|^p}\big)<\infty $ and $\sum_{k=2}^{\infty}( |a_{k}b_k|+2|a_{k}c_k|)<\infty.$ For example, if $p=2$ and let $a_k=\frac{1}{k^3},\,b_k=c_k=\frac{1}{k},\,k=2,3,\cdots.$
Then $\sum_{k=2}^{\infty}( |a_{k}b_k|+2|a_{k}c_k|)=\sum_{k=2}^{\infty}\frac{3}{k^2}<\infty.$
\end{example}
One can find similar examples such as the Hamilton-Jacobi equation in infinite dimensions (\cite{CL}) and the Laplacian defined on $\ell^2$ by Malliavin derivatives (\cite{Gro}).

Now, let us consider the Cauchy problem of the following first order linear homogenous partial differential equation  of infinitely many variables:
\begin{equation}\label{one order-zx}
\left\{
\begin{array}{ll}
 \displaystyle\partial_t u(t,\textbf{x}) -\sum_{i=1}^{\infty}a_{i}(t,\textbf{x}) \partial_{ x_i} u(t,\textbf{x}) -b(t,\textbf{x}) u(t,\textbf{x})= 0,\\[2mm]
 \displaystyle u(t,\textbf{x})\mid_{t=0}=\phi(\textbf{x}).
 \end{array}
\right.
 \end{equation}
Here $t\in \mathbb{R}, \textbf{x}=(x_i)_{i=1}^{\infty}\in \mathbb{R}^{\infty}$, the unknown $u$ is a real-valued function depending on $t$ and $\textbf{x}$, and the data $a_i(t,\textbf{x})$'s, $b(t,\textbf{x})$  and $\phi$ are analytic near $(0)$. Let
\begin{equation}\label{def of G}
G(t,\textbf{x}, \textbf{w})\triangleq\sum_{i=1}^{\infty}a_i(t,\textbf{x})w_i+w_0b(t,\textbf{x}),\quad \textbf{w}= (w_{ j})_{j=0}^{\infty}\in \mathbb{R}^{\infty}.
 \end{equation}
By modifying the proof of Theorem \ref{Infinite Dimensional Cauchy-Kowalevski Theorem}, we can show the following result.
\begin{corollary}\label{corollary 20}
Suppose $1\leq p<\infty$, the monomial expansion of $G $ near $(0)$ is absolutely convergent at a point near $\infty$ in the topology $\mathcal {T}^{p}$, and $D_\phi^{(0)}=\mathbb{R}^{\infty}$. Then there exists $r\in (0,\infty)$, independent of $\phi$, such that the equation (\ref{one order-zx}) admits locally an analytic solution $u$ near $(0)$ and $D_u^{(0)}\supset B_r^{p'}$. Furthermore, $u$ is Fr\'{e}chet differentiable with respect to $\ell^{p'}$ in $B_r^{p'}$, and the corresponding Fr\'{e}chet derivative $Du$ is continuous.
\end{corollary}

\section{Holmgren Type Theorem of Infinitely Many Variables}
Denote by $\Xi$ the set of real-valued functions, defined locally near (0) in $\mathbb{R}^{\infty}$, which are Fr\'{e}chet differentiable with respect to $\ell^2$ and whose Fr\'{e}chet derivative are locally continuous near $(0)$ in the $\mathscr{T}_2$ topology.

In this section, we shall establish the following Holmgren type theorem of infinitely many variables.


\begin{theorem}\label{Infinite Dimensional Holmgren Type Theorem}
Suppose the monomial expansion of $G$ (defined by (\ref{def of G})) near $(0)$ is absolutely convergent at a point near $\infty$ in the topology $\mathcal {T}^{\frac{1}{2}}$ and $\phi\in \Xi$. Then the solution to (\ref{one order-zx}) is locally unique in the class $\Xi$.
\end{theorem}

In the rest, we shall give a sketch of the proof of Theorem \ref{Infinite Dimensional Holmgren Type Theorem}.

It suffices to show that $U\in \Xi$ must be $0$ provided that $U$ solves (\ref{one order-zx}) with $\phi(\cdot)=0$. Without loss of generality, we assume that the data $a_i(t,\textbf{x})\equiv a_i(\textbf{x})$ ($i=1,2,\cdots$) and $b(t,\textbf{x})\equiv b(\textbf{x})$, i.e., they are independent of $t$. By our assumption, for each $i\in \mathbb{N}$, $a_{i}(\textbf{x})$ ({\it resp.} $b(\textbf{x})$) has a monomial expansion near $(0)$: $a_{i}(\textbf{x})=\sum_{\alpha\in \mathbb{N}^{(\mathbb{N})}}a_{i,\alpha}\textbf{x}^{\alpha}$ ({\it resp.} $b(\textbf{x})=\sum_{\alpha\in \mathbb{N}^{(\mathbb{N})}}a_{0,\alpha}\textbf{x}^{\alpha}$). Then, there exist $\rho_0>0$ and $0<s_j,t_i<1,\,0\leq i<\infty,\, 1\leq j<\infty$ such that $\sum_{j=1}^{\infty} s_j^{\frac{1}{2}}+\sum_{i=0}^{\infty}t_i^{\frac{1}{2}}  =1$ and
 \begin{equation}\label{equation 1133}
 \sum_{i=0}^{\infty} \left[\sum_{\alpha\in \mathbb{N}^{(\mathbb{N})}}|a_{i,\alpha}| (\frac{\rho_0}{s_j} )^{\alpha} \right]\frac{\rho_0}{t_i}<\infty.
  \end{equation}
Now we need to introduce a suitable transformation of variables. In a neighborhood of $(0)$ in $\mathscr{T}_2$, put
$$t^{'}\triangleq t+\sum_{i=1}^{\infty}  x_i^2,\quad x_i^{'}\triangleq x_i,\,\,i=1,2,\cdots, \quad\widetilde{U}(t^{'},x_1^{'},\cdots,x_i^{'},\cdots)\triangleq U\big(t^{'}-\sum_{i=1}^{\infty}  (x_i^{'})^2,x_1^{'},\cdots,x_i^{'},\cdots\big).
 $$
Then $\widetilde{U}\in\Xi$, and
\begin{eqnarray}\label{transformation 1}
 \bigg(1-2\sum_{i=1}^{\infty}x_i^{'}a_{i}(\textbf{x}^{'})\bigg)\partial_{t'} \widetilde{U}
 =\sum_{i=1}^{\infty}a_{i}(\textbf{x}^{'})\partial_{x_i^{'}} \widetilde{U}+b(\textbf{x}^{'})\widetilde{U},
\end{eqnarray}
where $\textbf{x}^{'}=(x_i^{'})\in\ell^2$. Clearly, $|2\sum_{i=1}^{\infty}x_i^{'}a_{i}(\textbf{x}^{'})|<1$ in a neighborhood of $(0)$ in $\mathcal {T}_2$. Thus (\ref{transformation 1}) can be written as $\partial_{t'} \widetilde{U}
 =\sum_{i=1}^{\infty}\widetilde{a}_{i}(\textbf{x}^{'})\partial_{x_i^{'}} \widetilde{U}+\widetilde{b}(\textbf{x}^{'})\widetilde{U},$ where $\widetilde{a}_{i}(\textbf{x}^{'})=\frac{a_{i}(\textbf{x}^{'})}{1-2\sum_{i=1}^{\infty}x_i^{'}a_{i}(\textbf{x}^{'})},1\,\,\leq i<\infty $ and $\widetilde{b}(\textbf{x}^{'})=\frac{b(\textbf{x}^{'})}{1-2\sum_{i=1}^{\infty}x_i^{'}a_{i}(\textbf{x}^{'})}.$ Now let $A_{0 } = t_0^{\frac{1}{4}},\,\,A_{i } =\max\Big\{t_i^{\frac{1}{4}},s_i^{\frac{1}{4}}\Big\},\,\,1\leq i<\infty.$ Then, $\sum_{i=0}^{\infty}A_{i}^2  < \sum_{j=1}^{\infty} s_j^{\frac{1}{2}}+\sum_{i=0}^{\infty}t_i^{\frac{1}{2}}  =1$ and by (\ref{equation 1133}) there exists some $\rho_1\in(0,\rho_0)$ such that $\sum_{i=0}^{\infty}\Big[\sum_{\alpha\in \mathbb{N}^{(\mathbb{N})}}|a_{i,\alpha}|\big(\frac{\rho_1}{A_j^4}\big)_{j\in \mathbb{N}}^{\alpha}\Big]\frac{\rho_1}{A_i^4}<\frac{1}{2}.$
Write the monomial expansions of $ \widetilde{b}(\textbf{x}^{'})$ and $ \widetilde{a}_{i}(\textbf{x}^{'}),\,1\leq i<\infty$ near $(0)$ respectively as $\widetilde{b}(\textbf{x}^{'})=\sum_{\alpha\in \mathbb{N}^{(\mathbb{N})}}a_{0,\alpha}^{'}\textbf{x}^{'\alpha}$ and $\widetilde{a}_{i}(\textbf{x}^{'})=\sum_{\alpha\in \mathbb{N}^{(\mathbb{N})}}a_{i,\alpha}^{'}\textbf{x}^{'\alpha},\,\,1\leq i<\infty.$
It is easily seen that
\begin{eqnarray}\label{convergence 1}
\sum_{i=0}^{\infty}\Bigg[\sum_{\alpha\in \mathbb{N}^{(\mathbb{N})}}|a_{i,\alpha}^{'}|\bigg(\frac{\rho_1}{A_j^4}\bigg)_{j\in \mathbb{N}}^{\alpha}\Bigg]\frac{\rho_1}{A_i^4}<\infty.
\end{eqnarray}
Formally write
$$
\begin{array}{ll}
\displaystyle G_1[\widetilde{U}]\triangleq\partial_{t'} \widetilde{U}
-\sum_{i=1}^{\infty}\widetilde{a}_{i}(\textbf{x}^{'})\partial_{x_i^{'}} \widetilde{U}-\widetilde{b}(\textbf{x}^{'})\widetilde{U},\\[2mm]
\displaystyle G_2[W]\triangleq-\partial_{t'} W
+\sum_{i=1}^{\infty}\partial_{x_i^{'}}[\widetilde{a}_{i}(\textbf{x}^{'})W]+ \Big [- \widetilde{b}(\textbf{x}^{'})
+ \frac{t^{'}}{A_0^3}-\sum_{i=1}^{\infty} \frac{x_i^{'}\widetilde{a}_{i}(\textbf{x}^{'})}{A_i^3}\Big]W,
\end{array}
$$
where the function $W$ will be defined later. For a sufficient small positive number $\lambda$, we put
$$H_{\lambda}= \{(t^{'},x_1^{'},\cdots,x_i^{'},\cdots)\in \mathbb{R}^{\infty}:\sum_{i=1}^{\infty}(x_i^{'})^2< t^{'}<\lambda \}.$$
Also, denote the boundary of $H_{\lambda}$ by $l_{\lambda}\cup I_{\lambda}\cup k_{\lambda}$ where
$$
\begin{array}{ll}
\displaystyle l_{\lambda}=\{(t^{'},x_1^{'},\cdots,x_i^{'},\cdots)\in \mathbb{R}^{\infty}: t^{'}=\lambda, \sum_{i=1}^{\infty}(x_i^{'})^2<\lambda\},\\[2mm]
\displaystyle  I_{\lambda}=\{(t^{'},x_1^{'},\cdots,x_i^{'},\cdots)\in \mathbb{R}^{\infty}:  t^{'}= \lambda,\sum_{i=1}^{\infty}(x_i^{'})^2= t^{'}\},\\[2mm]
\displaystyle k_{\lambda}=\{(t^{'},x_1^{'},\cdots,x_i^{'},\cdots)\in \mathbb{R}^{\infty}:  t^{'}< \lambda,\sum_{i=1}^{\infty}(x_i^{'})^2= t^{'}\}.
\end{array}
$$
From the definition one can see that $H_{\lambda}$ is an open subset of $\ell^2$ and $\widetilde{U}\mid_{I_{\lambda}\cup k_{\lambda}}=0$. Let
\begin{eqnarray}\label{def of F}
 F(t^{'},\textbf{x}^{'}) \triangleq\bigg(\frac{t^{'}}{A_0},-\frac{\widetilde{a}_{1}(\textbf{x}^{'}) }{A_1 },\cdots,-\frac{\widetilde{a}_{i}(\textbf{x}^{'}) }{A_i },\cdots\bigg), \qquad t^{'}\in\mathbb{R}, \;\textbf{x}^{'}=(x_i^{'})\in\ell^2,
\end{eqnarray}
$H\triangleq\Big\{(t^{'},x_1^{'},\cdots,x_i^{'},\cdots)\in \mathbb{R}^{\infty}: \big(\frac{t^{'}}{A_0 }\big)^2+ \sum_{i=1}^{\infty}\big(\frac{x_i^{'}}{A_i }\big)^2<\infty\Big\}$ and $B\triangleq \{(t^{'},x_1^{'},\cdots,x_i^{'},\cdots)\in \mathbb{R}^{\infty}:  ( t^{'} )^2+ \sum_{i=1}^{\infty} ( x_i^{'} )^2<\infty \}.$ We can view $B$ as $\mathbb{R}\times \ell^2$ and we will use this convention later. Since $\sum_{i=0}^{\infty} A_i^2<\infty,$ we have the natural inclusion map $i$ which maps $H$ into $B$ such that the triple $(i,H,B)$ is an abstract Wiener space. Since $H=H^*$, $i^{*}B^{*}$ can be characterized by a subset of $H$ which will be denoted by $B^{*}$. Precisely, $B^{*}\triangleq \{h\in H:  \,\text{ there exists }C_h\in (0,+\infty)\text{ such that }|(h,g)_{H}|\leq C_h||ig||_{B}\text{ for any }g\in H \}.$ A simple computation shows that $B^{*}= \Big\{(t^{'},x_1^{'},\cdots,x_i^{'},\cdots):  \frac{(t^{'})^2}{A_0^4 }+ \sum_{i=1}^{\infty} \frac{(x_i^{'})^2}{A_i^4 } <\infty \Big\}.$ In order to apply Theorem \ref{Divergence Theorem}, we need the following result:
\begin{proposition}\label{3 cond F satisfy}
\begin{itemize}For a sufficiently small $\lambda>0$ it holds that
\item[(a)]
$F$ is uniformly bounded in $B^*$ norm, i.e., $\sup_{(t^{'},\textbf{x}^{'})\in H_{\lambda}}\sum_{i=1}^{\infty}  \Big[\frac{\widetilde{a}_{i}(\textbf{x}^{'}) }{A_i^3} \Big]^2 <\infty;$
\item[(b)] $F$ maps $H_{\lambda}\cup l_{\lambda}\cup k_{\lambda}\cup I_{\lambda}$ into a bounded subset of $H$, i.e., $\sup_{(t^{'},\textbf{x}^{'})\in H_{\lambda}\cup l_{\lambda}\cup k_{\lambda}\cup I_{\lambda}}\sum_{i=1}^{\infty} \Big[\frac{\widetilde{a}_{i}(\textbf{x}^{'}) }{A_i^2} \Big]^2 <\infty;$
\item[(c)]
$F$ is $H$-Fr\'{e}chet differentiable with $DF$ being a trace class operator and the trace norm $||DF||_1$ is $p$ integrable on $ H_{\lambda}$.
\end{itemize}
\end{proposition}

It is easy to see that both $l_{\lambda}$ and $k_{\lambda}$ are two differentiable surfaces in $B$, by \cite[Remark 2]{Goo}, both $l_{\lambda}$ and $k_{\lambda}$ are two $H$-$C^1$ surfaces in the abstract Wiener space $(i,H,B)$. Denote by $p$ ({\it resp.} $\sigma$) the corresponding Wiener measure on $B$ ({\it resp.} normal surface measure given in Theorem \ref{theorem 1 in Goodman DT}) with parameters $x=0$ and $t=1$. We may show the following result (which means that, for any $\lambda>0$,  $\sigma$ is a finite measure):
\begin{lemma}\label{small surface measure}
For any $\lambda>0$, it holds that $\sigma(l_{\lambda}\cup k_{\lambda}\cup I_{\lambda})<\infty$.
\end{lemma}

By Lemma \ref{small surface measure} and the second conclusion in  Proposition \ref{3 cond F satisfy}, it follows that $|F|_H$ is $\sigma$-integrable on $l_{\lambda}\cup k_{\lambda}\cup I_{\lambda}$.

Note that the Divergence Theorem, i.e., Theorem \ref{Divergence Theorem} requires the boundary to be a ``smooth" surface but the boundary of $H_{\lambda}$ is union of two ``smooth" surfaces. However, similar arguments in \cite{Goo} can be modified to show that the Divergence Theorem also holds in this following case:

\begin{theorem}\label{Divergence theorem coro}
Assume that ${\cal F}: H_{\lambda}\cup l_{\lambda}\cup k_{\lambda}\cup I_{\lambda}\rightarrow H$ is a continuous function with finite divergence on $H_{\lambda}$ and that ${\cal F}$ is uniformly bounded with respect to the $B^*$-norm on $H_{\lambda}$. If the function $|{\cal F}(\cdot)|$ is integrable with respect to the normal surface measure $\sigma( \cdot)$ on $l_{\lambda}\cup k_{\lambda}\cup I_{\lambda}$, and the trace class norm of $D{\cal F}$ is $p$-integrable on $H_{\lambda}$, then
\begin{eqnarray*}
\int_{H_{\lambda}}\left[({\rm div }\;{\cal F})(y)-\langle {\cal F}(y),y\rangle\right]\,p(\mathrm{d}y)=\int_{l_{\lambda}\cup k_{\lambda}\cup I_{\lambda}}\langle {\cal F}(y),{\bf n}(y)\rangle\,\sigma(\mathrm{d}y).
\end{eqnarray*}
\end{theorem}

By Proposition \ref{3 cond F satisfy} and  Lemma \ref{small surface measure}, we conclude that $H_{\lambda}$ and the function $F$ defined by (\ref{def of F}) satisfies the assumptions in Theorem \ref{Divergence theorem coro}.  
If $W $ is Fr\`{e}chet differentiable respect to $B$ with continuous Fr\`{e}chet derivatives, then $F_0=W \widetilde{U} F$ satisfies the assumptions in Theorem \ref{Divergence theorem coro}.  Clearly, $\text{div}(W \widetilde{U} F)=A_0\partial_{t'}  \Big[\frac{  W\widetilde{U}}{A_0} \Big]
-\sum_{i=1}^{\infty}A_i\partial_{x_i^{'}}\Big[\frac{\widetilde{a}_{i} W\widetilde{U}}{A_i} \Big]$ and
$\langle W \widetilde{U} F(y),y \rangle= \frac{W\widetilde{U}\cdot t^{'}}{A_0^3}-\sum_{i=1}^{\infty}\frac{\widetilde{a}_{i}(\textbf{x}^{'})W\widetilde{U}\cdot x_i^{'}}{A_i^3},$ where $y=(t^{'}, \textbf{x}^{'})=(t^{'},(x_i^{'})_{i=1}^{\infty} )\in H$. Apply Theorem \ref{Divergence theorem coro}, we have
\begin{eqnarray}
 \int_{H_{\lambda}}(WG_1[\widetilde{U}]- \widetilde{U}G_2[W])\,p(\mathrm{d}y)
&=&\frac{\lambda}{A_0^2} \int_{l_{\lambda}}W\widetilde{U} \,\sigma(\mathrm{d}y)\label{formula 20}.
\end{eqnarray}
Let $b'(\textbf{x}^{'})\triangleq \Big[ \sum_{i=1}^{\infty}\partial_{x_i^{'}}\widetilde{a}_{i}(\textbf{x}^{'})- \widetilde{b}(\textbf{x}^{'})
+ \frac{t^{'}}{A_0^3}-\sum_{i=1}^{\infty} \frac{x_i^{'}\widetilde{a}_{i}(\textbf{x}^{'})}{A_i^3}\Big].$ Then $G_2[W]=-\partial_{t'} W
+\sum_{i=1}^{\infty}\widetilde{a}_{i}(\textbf{x}^{'})\partial_{x_i^{'}}W+b'(\textbf{x}^{'})W.$ Let $G' (\textbf{x}', \textbf{w} )\triangleq\sum_{i=1}^{\infty}\widetilde{a}_{i}(\textbf{x}')w_i+w_0b(\textbf{x}'),$ where $\textbf{w}= (w_{ j})_{j=0}^{\infty}\in \mathbb{R}^{\infty}$ and $\textbf{x}= (x_i)_{i=1}^{\infty}\in \mathbb{R}^{\infty}.$ Recall that
$\sum_{i=0}^{\infty}\Big[\sum_{\alpha\in \mathbb{N}^{(\mathbb{N})}}|a_{i,\alpha}^{'}|\big(\frac{\rho_1}{A_j^4}\big)^{\alpha}\Big]\frac{\rho_1}{A_i^4}<\infty.$ Hence, for any $\rho\in (0,\rho_1)$, the monomial expansion of $G'$ is absolutely convergent at $\big(\big(\frac{\rho}{A_j}\big)_{j=0}^{\infty},\big(\frac{\rho}{A_i}\big)_{i=1}^{\infty}\big),$ which is a point near $\infty$ in the topology $\mathcal {T}^{2}$ by the fact that $\sum_{i=0}^{\infty}A_{i}^2  < 1.$ Therefore, the equation $G_2[W]=0$ satisfies the assumptions in Corollary \ref{corollary 20}. Hence, for any $n\geq 0$ and $k_1,\cdots,k_n\in \mathbb{N}$, there is an analytic solution $W$ to $G_2[W]=0$ and $W\mid_{t^{'}=\lambda}=(x_1^{'})^{k_1}\cdots(x_n^{'})^{k_n}$. From Corollary \ref{corollary 20}, it follows that the monomial expansion of $W$ is absolutely convergent in a neighborhood of $(0)$ containing $H_{\lambda}$ for sufficiently small $\lambda$ and $W$ is Fr\`{e}chet differentiable with respect to $B$ and the Fr\`{e}chet derivative is continuous in this neighborhood. Therefore, for any sufficiently small $\lambda>0$, applying (\ref{formula 20}), we arrive at
\begin{eqnarray}\label{polynomial is dense}
\int_{l_{\lambda}}(x_1^{'})^{k_1}\cdots(x_n^{'})^{k_n}\widetilde{U} \,\sigma(\mathrm{d}y)=0,
\end{eqnarray}
for any $n\geq 0$ and $k_1,\cdots,k_n\in \mathbb{N}$. Let $L_{\lambda}= \{(t^{'},x_1^{'},\cdots,x_i^{'},\cdots)\in \mathbb{R}^{\infty}: t^{'}=\lambda, \sum_{i=1}^{\infty}(x_i^{'})^2<\infty \},$ then the surface measure $\sigma^{'}$ on $L_{\lambda}$ is identified with the Gaussian measure with parameters $x=0$ and $t=1$ on the Hilbert space $H_0=\{(x_1^{'},\cdots,x_i^{'},\cdots)\in \mathbb{R}^{\infty}:  \sum_{i=1}^{\infty}(x_i^{'})^2<\infty \}.$ Note that (\ref{polynomial is dense}) is equivalent to
\begin{eqnarray*}
\int_{L_{\lambda}}(x_1^{'})^{k_1}\cdots(x_n^{'})^{k_n}\chi_{l_{\lambda}}\widetilde{U} \,\sigma^{'}(\mathrm{d}y)=0.
\end{eqnarray*}
We also need the following density result.
\begin{lemma}\label{dense lemma}
$\span\{(x_1^{'})^{k_1}\cdots(x_n^{'})^{k_n}:n\geq 0,\, k_1,\cdots,k_n\in \mathbb{N}\}$ is dense in $L^2(L_{\lambda},\sigma^{'})$.
\end{lemma}

One may check that $L^2(l_{\lambda},\sigma)=\chi_{l_{\lambda}}\cdot L^2(L_{\lambda},\sigma^{'})$.
From Lemma \ref{dense lemma} and noting the continuity of $\widetilde{U}$, we deduce that $\widetilde{U}\equiv 0$ on $l_{\lambda}$ for any sufficiently small $\lambda>0$ and hence $\widetilde{U}\mid_{k_{\lambda}}=0$. Therefore, $\widetilde{U}\equiv 0$ on $H_{\lambda}$ for any sufficiently small $\lambda>0$. This implies that $U\equiv 0$ at a neighborhood of $(0)$ in the $\mathcal {T}_2$ topology restricted on the half space $t>0$. By the same way we can prove that $U\equiv 0$ at a neighborhood of $(0)$ the $\mathcal {T}_2$ topology restricted on the half space $t<0$. Finally, by the continuity of $U$ we have $U\equiv 0$ at a neighborhood of $(0)$ in the $\mathcal {T}_2$ topology.


\Acknowledgements{ The research of the first author is supported by NSFC  under grant 11501384;
The research of the second author is supported by NSFC under grant 11221101, the NSFC-CNRS
Joint Research Project under grant 11711530142 and the PCSIRT under
grant IRT$\_$16R53 from the Chinese Education
Ministry.}

\end{document}